\documentclass[smallcondensed]{svjour3}

\smartqed

\usepackage{amsmath}
\usepackage{amssymb}

\usepackage{graphicx}

\newcommand{\I}{\mathbf I}
\newcommand{\M}{\mathbf M}
\newcommand{\N}{\mathbf N}

\newcommand{\T}{\mathbf T}
\newcommand{\Tss}{\mathbf{T}_{ss}}
\newcommand{\pp}{\mathbf{e}_1}
\newcommand{\qq}{\mathbf{e}_2}

\newcommand{\Xt}{\mathbf X_t}
\newcommand{\X}{\mathbf X}
\newcommand{\Xs}{\mathbf X_s}
\newcommand{\Xss}{\mathbf X_{ss}}
\newcommand{\Tt}{\mathbf T_t}
\newcommand{\jacobi}[2]{\left(\frac{#1}{#2}\right)}

\journalname{Acta Appl. Math.}

\begin{document}

\title{The vortex filament equation as a pseudorandom generator\thanks{This work was supported by the Spanish Ministry of Economy and Competitiveness, with the project MTM2011-24054, and by the Basque Government, with the project IT641-13.}}

\author{Francisco de la Hoz \and Luis Vega}

\institute{Francisco de la Hoz \at
              Department of Applied Mathematics and Statistics and Operations Research, Faculty of Science and Technology, University of the Basque Country UPV/EHU, Barrio Sarriena S/N, 48940 Leioa, Spain \\
              \email{francisco.delahoz@ehu.es}
           \and
           Luis Vega \at
              Department of Mathematics, Faculty of Science and Technology, University of the Basque Country UPV/EHU, Barrio Sarriena S/N, 48940 Leioa, Spain \\ BCAM - Basque Center for Applied Mathematics, Alameda Mazarredo, 14, 48009 Bilbao, Spain \\ \email{luis.vega@ehu.es; lvega@bcamath.org}
}

\date{Received: date / Accepted: date}

\maketitle

\begin{abstract}

\noindent In this paper, we consider the evolution of the so-called vortex filament equation (VFE),
\begin{equation*}
\mathbf X_t = \mathbf X_s\wedge\mathbf X_{ss},
\end{equation*}

\noindent taking a planar regular polygon of $M$ sides as initial datum. We study VFE from a completely novel point of view: that of an evolution equation which yields a very good generator of pseudorandom numbers in a completely natural way. This essential randomness of VFE is in agreement with the randomness of the physical phenomena upon which it is based.



\end{abstract}

\section{Introduction}

The binormal flow,
\begin{equation}
\label{e:binormal}
\Xt = \kappa\mathbf b,
\end{equation}

\noindent where $t$ is the time, $\kappa$ the curvature, and $\mathbf b$ the binormal component of the
Frenet-Serret formulae, appeared for the first time in 1906 \cite{darios}, and was rederived in \cite{arms}, as an approximation of the dynamics of a vortex filament under the Euler equations. It is also known as the vortex filament equation (VFE) or the localized induction equation (LIA). An equivalent expression of \eqref{e:binormal} is
\begin{equation}
\label{e:xt}\Xt = \Xs\wedge\Xss,
\end{equation}

\noindent where $\wedge$ is the usual cross-product, and $s$ is the arc-length parameter. The tangent vector $\T = \X_s$ remains with constant length and, hence, we can assume that $\|\T\|_2 = 1$, for all time. Differentiating \eqref{e:xt} with respect to $s$, we get
\begin{equation}
\label{e:schmap}\Tt = \T\wedge\Tss,
\end{equation}

\noindent known as the Schr\"odinger map on the sphere.

The question of making sense of initial data with corners in \eqref{e:xt}-\eqref{e:schmap} has recently received some attention. For instance, the existence of solutions starting with a single corner, which are precisely the self-similar solutions of \eqref{e:xt}-\eqref{e:schmap}, has been proven in \cite{gutierrez} (see also \cite{delahoz2007} for the corresponding problem in the hyperbolic space); and numerical simulations of these solutions have been carried out in \cite{buttke87,DelahozGarciaCerveraVega09}). Furthermore, the fact that this kind of solutions yields a well-posed problem has been shown in a long-term collaboration between Banica and Vega \cite{BV0,BV1,BV2,BV3}; in particular, the last paper of the series, \cite{BV3}, closes the question, because it proves that the problem with single-corner initial data is well-posed in an adequate function space.

Even if the solutions of \eqref{e:xt}-\eqref{e:schmap} for single-corner initial data are well understood, very little has been done for more general initial data with several corners \cite{Didier}. However, in a recently submitted paper \cite{HV2013}, we have studied for the first time the evolution of \eqref{e:xt}-\eqref{e:schmap}, taking a regular planar polygon of $M$ sides as the initial datum. The main ideas of \cite{HV2013} are as follows. In order to avoid working with the curvature $\kappa$ and the torsion $\tau$, we consider an alternative version of the Frenet-Serret formulae,
\begin{equation}
\label{e:Te1e2}
\begin{pmatrix}
\T \cr \pp \cr \qq
\end{pmatrix}_s =
\begin{pmatrix}
    0 & \alpha & \beta
    \cr
    - \alpha & 0 & 0 \cr -\beta & 0
& 0
\end{pmatrix} \cdot
\begin{pmatrix}
    \T \cr \pp \cr \qq
\end{pmatrix};
\end{equation}

\noindent where
\begin{equation}
\label{e:abkt}
\alpha(s,t) = \kappa(s, t)\cos\left(\int^s\tau(s',t)ds'\right), \quad \beta(s,t) = \kappa(s, t)\sin\left(\int^s\tau(s',t)ds'\right).
\end{equation}

\noindent Then, the Hasimoto transformation \cite{hasimoto} adopts the form
\begin{equation}
\psi = \alpha + i\beta,
\end{equation}

\noindent and transforms \eqref{e:xt}-\eqref{e:schmap} into the nonlinear Schr\"odinger (NLS) equation:
\begin{equation}
\label{e:schr}
\psi_t = i\psi_{ss} + i\left(\frac{1}{2}(|\psi|^2 + A(t))\right)\psi,
\end{equation}

\noindent where $A(t)$ is a certain time-dependent real constant. The main idea is to work with \eqref{e:schr}, and, at a given $t$, to recover $\X(s, t)$ and $\T(s, t)$ from $\psi(s, t)$ by integrating \eqref{e:Te1e2}, up to a rigid movement that can be determined by the symmetries of the problem.

Observe that, if we define $\N\equiv \pp + i\qq$ (see \cite{hasimoto}), then it is not difficult to check that \eqref{e:schmap} can be rewritten as
\begin{equation}
\label{e:TtN}
\Tt = \frac{i}{2}(\psi_s\bar\N - \bar\psi_s\N).
\end{equation}

\noindent Therefore, if the system $\{\psi,\T,\N\}$ solves \eqref{e:TtN}, then, defining
\begin{equation}
\label{e:hasimoto2}
\tilde\psi(s, t) \equiv e^{i\omega(t)}\psi(s, t),
\end{equation}

\noindent $\{\tilde\psi,\T, e^{i\omega(t)}\N\}$ is also a solution; i.e., the tangent vector $\T$ does not change, while the vectors $\pp$ and $\qq$ rotate in the normal plane $\omega(t)$ degrees around $\T$. Since in this paper (and in \cite{HV2013}) we are interested only in $\T$, we conclude that $\psi(s, t)$ can be chosen without loss of generality up to a complex value (that depends on time) with modulus one; in particular, we can choose $\psi(s, t)$ to be real, $\psi(s, t)\equiv|\psi(s, t)|$.

Given a regular planar polygon of $M$ sides as $\X(s, 0)$, there is no torsion; hence, from \eqref{e:abkt}, $\psi(s, 0)$ is precisely the curvature of the polygon, which is a $2\pi/M$-periodic sum of Dirac deltas:
\begin{equation}
\label{e:curvatureM}
\psi(s, 0) \equiv \kappa(s) = \frac{2\pi}{M}\sum_{k = -\infty}^{\infty}\delta(s - \tfrac{2\pi k}{M}).
\end{equation}

\noindent Then, bearing in mind the Galilean invariance of \eqref{e:schr} and \textbf{assuming uniqueness}, we are able to obtain $\psi(s, t)$ at any rational multiple of $2\pi / M^2$. During all this paper, we assume that $p$ and $q$ are two \textbf{coprime} natural numbers. Defining $t_{pq} \equiv (2\pi/M^2)(p/q)$, it can be shown that
\begin{equation}
\label{e:psistpq}
\psi(s, t_{pq}) = \frac{2\pi}{Mq}\hat\psi(0, t_{pq})\sum_{k = -\infty}^\infty\sum_{m = 0}^{q - 1}G(-p, m, q) \delta(s - \tfrac{2\pi k}{M} - \tfrac{2\pi m}{Mq}),
\end{equation}

\noindent where $\hat\psi(0, t_{pq})$ is the mean of $\psi(s, t_{pq})$ over a period,
\begin{equation}
\hat\psi(0, t_{pq}) = \frac{M}{2\pi}\int_0^{2\pi/M}\psi(s, t_{pq})ds,
\end{equation}

\noindent and
\begin{equation}
G(a, b, c) = \sum_{l=0}^{c - 1}e^{2\pi i (al^2 + bl)/c}
\end{equation}

\noindent denotes a generalized quadratic Gau{\ss} sum. Remark that, as explained in the lines following \eqref{e:hasimoto2}, we can assume without loss of generality that $\hat\psi(0, t_{pq})$ is real.

An important property of the generalized quadratic Gau{\ss} sums is that
\begin{equation}
\label{e:G-pmq}
|G(-p, m, q)| =
\begin{cases}
\sqrt q, & \mbox{if $q\equiv 1\bmod2$},
    \\
\sqrt{2q}, & \mbox{if $q\equiv 0\bmod2$ $\wedge$ $q/2\equiv m\bmod 2$},
    \\
0, & \mbox{if $q\equiv 0\bmod2$ $\wedge$ $q/2\not\equiv m\bmod 2$};
\end{cases}
\end{equation}

\noindent therefore, we can write
\begin{equation}
\label{e:G-pmq2}
G(-p, m, q) =
\begin{cases}
\sqrt q e^{i \theta_m}, & \mbox{if $q\equiv 1\bmod2$},
    \\
\sqrt{2q} e^{i \theta_m}, & \mbox{if $q\equiv 0\bmod2$ $\wedge$ $q/2\equiv m\bmod 2$},
    \\
0, & \mbox{if $q\equiv 0\bmod2$ $\wedge$ $q/2\not\equiv m\bmod 2$},
\end{cases}
\end{equation}

\noindent for certain $\theta_m$ that also depend on $q$. Hence, defining
\begin{equation}
\label{e:rho}
\rho =
\begin{cases}
\frac{2\pi}{M\sqrt q}\hat\psi(0, t_{pq}), & \mbox{if $q\equiv1\bmod2$},
    \cr
\frac{2\pi}{M\sqrt {\tfrac{q}{2}}}\hat\psi(0, t_{pq}), & \mbox{if $q\equiv0\bmod2$ $\wedge$ $q/2\equiv m\bmod 2$},
    \cr
0, & \mbox{if $q\equiv0\bmod2$ $\wedge$ $q/2\not\equiv m\bmod 2$},
\end{cases}
\end{equation}

\noindent we represent \eqref{e:psistpq} as
\begin{equation}
\label{psistpq}
\psi(s, t_{pq}) = \sum_{k=-\infty}^{\infty}\sum_{m = 0}^{q - 1}\rho e^{i \theta_m} \delta(s - \tfrac{2\pi k}{M} - \tfrac{2\pi m}{Mq}).
\end{equation}

\noindent The coefficients multiplying the Dirac deltas are in general not real, except for $t = 0$ and $t_{1,2} = \pi/M^2$. Therefore, $\psi(s, t_{pq})$ does not correspond to a planar polygon, but to a skew polygon with $Mq$ (for $q$ odd) or $Mq/2$ (for $q$ even) equal-lengthed sides.

In order to recover $\X$ and $\T$ from $\psi$, we observe that every addend $\rho e^{i\theta_m}\delta(s - \tfrac{2\pi m}{Mq})$ in \eqref{psistpq}, with $\rho\not=0$, induces a rotation on $\T$, $\pp$ and $\qq$. More precisely, defining $c_\rho \equiv \cos(\rho)$, $s_\rho \equiv \sin(\rho)$, $c_{\theta_m} \equiv \cos(\theta_m)$, $s_{\theta_m} \equiv \sin(\theta_m)$,
\begin{equation}
\label{e:Mm}
\M_m =
\begin{pmatrix}
c_\rho & s_\rho c_{\theta_m} & s_\rho s_{\theta_m}
    \\
-s_\rho c_{\theta_m} & c_\rho c_{\theta_m}^2 + s_{\theta_m}^2 & (c_\rho - 1)c_{\theta_m}s_{\theta_m}
    \\
-s_\rho s_{\theta_m} & (c_\rho - 1)  c_{\theta_m}s_{\theta_m} & c_\rho s_{\theta_m}^2 + c_{\theta_m}^2
\end{pmatrix}
\end{equation}

\noindent is the matrix such that
\begin{equation}
\label{e:T2pik/mq}
\left(
    \begin{array}{c}
\T(\tfrac{2\pi m}{Mq}^+)
    \cr
\hline
\pp(\tfrac{2\pi m}{Mq}^+)
    \cr
\hline
\qq(\tfrac{2\pi m}{Mq}^+)
    \end{array}
\right)
    = \M_m\cdot
\left(
    \begin{array}{c}
\T(\tfrac{2\pi m}{Mq}^-)
    \cr
\hline
\pp(\tfrac{2\pi m}{Mq}^-)
    \cr
\hline
\qq(\tfrac{2\pi m}{Mq}^-)
    \end{array}
\right),
\end{equation}

\noindent where all the vectors are row vectors. Notice that, when $\rho = 0$, $\M_m$ is just the identity matrix $\I$. From \eqref{e:Mm}, it follows that the non-zero value of $\rho$ is the angle between any two adjacent sides. Imposing that \eqref{e:psistpq} corresponds to a closed polygon, i.e., that
\begin{equation}
\label{e:productM}
\M_{Mq-1}\cdot \M_{Mq-2} \cdot \ldots \cdot \M_1 \cdot \M_0 \equiv \I,
\end{equation}

\noindent there is very strong evidence that the non-zero value of $\rho$ is given by
\begin{equation}
\label{e:cosrho}
\cos(\rho) =
\begin{cases}
2\cos^{2/q}(\tfrac{\pi}{M}) - 1, & \mbox{if $q\equiv1\bmod2$},
    \cr
2\cos^{4/q}(\tfrac{\pi}{M}) - 1, & \mbox{if $q\equiv0\bmod2$};
\end{cases}
\end{equation}

\noindent and the value of $\hat\psi(0, t_{pq})$ follows from \eqref{e:rho}.

The previous ideas suggest very strongly that $\psi(s, t)$ is also periodic in time, with period $2\pi / M^2$. Furthermore, bearing in mind the symmetries of the problem, it follows that also $\T$ is periodic in time, while $\X$ is periodic in time up to a movement of its center of mass with constant upward velocity.

Although the study of VFE is interesting per se, a recurring question is up to what extent it is valid as a simplified model for describing real vortex filament motion. In this paper, we would like to make a step forward in that direction, by proving that the evolution of $\X$ and $\T$ for a regular polygonal initial datum is essentially random, or, in other words, that it gives as a by-product a simple and powerful generator of pseudorandom numbers. More precisely, fixed $q$, we will focus on two quantities: the triple product of three consecutive tangent vectors, and the scalar product of a tangent vector and the second next one. Furthermore, taking these two quantities respectively as the real and imaginary parts of a complex number, we will have a generator of pseudorandom numbers located on a circumference of center $ic_\rho^2$ and radius $s_\rho^2$.

The structure of this paper is as follows. In Section \ref{s:quantities}, we study the aforementioned quantities. We prove that they depend exclusively on $\phi(p)$, which is defined as the inverse of a multiple of $p$ in a finite ring:
\begin{equation}
\phi(p) \equiv
\begin{cases}
(4p)^{-1}\bmod q, & \mbox{if $q\equiv1\bmod 2$},
    \\
p^{-1}\bmod (q/2), & \mbox{if $q\equiv2\bmod 4$},
    \\
p^{-1}\bmod q, & \mbox{if $q\equiv0\bmod 4$}.
\end{cases}
\end{equation}

\noindent Therefore, it is convenient to consider three cases of growing difficulty, according to the oddness of $q$ and $q/2$: Section \ref{s:q12} deals with $q$ odd; Section \ref{s:q24} deals with $q$ even, but $q/2$ odd; and Section \ref{s:q04} deals with both $q$ and $q/2$ even.

In Section \ref{s:randomness}, we analyze the pseudorandom properties of $\phi(p)$. More precisely, we put it in the frame of the so-called explicit inversive congruential generators. Finally, in Section \ref{s:conclusions}, we draw the main conclusions and point out future directions to extend this research.

\section{Two interesting quantities}

\label{s:quantities}

As we have mention in the introduction, we will divide the problem in three cases, according to the oddness of $q$ and $q/2$.

\subsection{Case with $q\equiv1\bmod2$}

\label{s:q12}

The simplest case is when $q$ is odd. Then, $\psi(s, t_{pq})$ in \eqref{psistpq} adopts over the first period the form
\begin{equation}
\psi(s, t_{pq}) = \rho\sum_{m = 0}^{q - 1}e^{i \theta_m} \delta(s - \tfrac{2\pi m}{Mq}), \quad s\in[0, \tfrac{2\pi}{M}),
\end{equation}

\noindent i.e., the vertices of $\X$, denoted by $\X_m$, are located at $s = \tfrac{2\pi m}{Mq}$, and the sides are the segments that join $\X_{m + 1}$ and $\X_m$. As stated in the introduction, we are interested in calculating the triple product of $\T(\tfrac{2\pi m}{Mq}^-)$, $\T(\tfrac{2\pi m}{Mq}^+)\equiv\T(\tfrac{2\pi (m + 1)}{Mq}^-)$, and $\T(\tfrac{2\pi (m+1)}{Mq}^+)$; and the scalar product of $\T(\tfrac{2\pi m}{Mq}^-)$ and $\T(\tfrac{2\pi (m+1)}{Mq}^+)$. Let us calculate the first quantity:
\begin{align}
\label{e:phirho}
\left[\T(\tfrac{2\pi m}{Mq}^-), \T(\tfrac{2\pi m}{Mq}^+), \T(\tfrac{2\pi (m+1)}{Mq}^+) \right]
& =
\left(\T(\tfrac{2\pi m}{Mq}^-)\wedge \T(\tfrac{2\pi m}{Mq}^+)\right)\cdot\T(\tfrac{2\pi (m+1)}{Mq}^+)
    \cr
& =
\begin{vmatrix}
\T(\tfrac{2\pi m}{Mq}^-)
    \cr
\T(\tfrac{2\pi m}{Mq}^+)
    \cr
\T(\tfrac{2\pi (m+1)}{Mq}^+)
\end{vmatrix}.
\end{align}

\noindent It is important to bear in mind that both the triple product of three vectors and the scalar product of two vectors are rotation-invariant, and, thus, we do not have to determine the global rotation of the whole skew polygon, which is very involved. Instead, we can simply assume that $\T(\tfrac{2\pi m}{Mq}^-) = (1,0,0)$, $\pp(\tfrac{2\pi m}{Mq}^-) = (0,1,0)$, and $\qq(\tfrac{2\pi m}{Mq}^-) = (0,0,1)$, i.e., they form the identity matrix. Then, from \eqref{e:Mm}-\eqref{e:T2pik/mq},
\begin{equation}
\left(
    \begin{array}{c}
\T(\tfrac{2\pi (m + 1)}{Mq}^-)
    \cr
\hline
\pp(\tfrac{2\pi (m + 1)}{Mq}^-)
    \cr
\hline
\qq(\tfrac{2\pi (m + 1)}{Mq}^-)
    \end{array}
\right)
    =
\left(
    \begin{array}{c}
\T(\tfrac{2\pi m}{Mq}^+)
    \cr
\hline
\pp(\tfrac{2\pi m}{Mq}^+)
    \cr
\hline
\qq(\tfrac{2\pi m}{Mq}^+)
    \end{array}
\right)
    = \M_m,
\end{equation}

\noindent and
\begin{equation}
\left(
    \begin{array}{c}
\T(\tfrac{2\pi (m + 1)}{Mq}^+)
    \cr
\hline
\pp(\tfrac{2\pi (m + 1)}{Mq}^+)
    \cr
\hline
\qq(\tfrac{2\pi (m + 1)}{Mq}^+)
    \end{array}
\right)
    = \M_{m + 1} \cdot
\left(
    \begin{array}{c}
\T(\tfrac{2\pi (m + 1)}{Mq}^-)
    \cr
\hline
\pp(\tfrac{2\pi (m + 1)}{Mq}^-)
    \cr
\hline
\qq(\tfrac{2\pi (m + 1)}{Mq}^-)
    \end{array}
\right).
\end{equation}

\noindent More precisely, $\T(\tfrac{2\pi m}{Mq}^+)$ is the first row of $\M_m$, while $\T(\tfrac{2\pi (m+1)}{Mq}^+)$ is the first row of $\M_{m + 1}\cdot \M_m$. Defining $\Delta_m = \theta_{m + 1} - \theta_{m}$, $c_{\Delta_m} = \cos(\Delta_m)$, $s_{\Delta_m} = \sin(\Delta_m)$, a straight calculation shows that
\begin{equation}
\label{e:Tm+1}
\T(\tfrac{2\pi (m+1)}{Mq}^+) =
\begin{pmatrix}
c_\rho^2 -s_\rho^2c_{\Delta_m}
    &
\begin{array}{r}c_\rho s_\rho c_{\theta_m}(1 + c_{\Delta_m}) \\ - s_\rho s_{\theta_m}s_{\Delta_m}\end{array}
    &
\begin{array}{r}c_\rho s_\rho s_{\theta_m}(1 + c_{\Delta_m}) \\ + s_\rho c_{\theta_m}s_{\Delta_m}\end{array}
\end{pmatrix}.
\end{equation}

\noindent Therefore, \eqref{e:phirho} becomes
\begin{align}
\begin{vmatrix}
\T(\tfrac{2\pi m}{Mq}^-)
    \cr
\T(\tfrac{2\pi m}{Mq}^+)
    \cr
\T(\tfrac{2\pi (m+1)}{Mq}^+)
\end{vmatrix}
& =
\begin{vmatrix}
1 & 0 & 0
    \\
c_\rho & s_\rho c_{\theta_m} & s_\rho s_{\theta_m}
    \\
c_\rho^2 -s_\rho^2c_{\Delta_m}
    &
\begin{array}{r}c_\rho s_\rho c_{\theta_m}(1 + c_{\Delta_m}) \\ - s_\rho s_{\theta_m}s_{\Delta_m}\end{array}
    &
\begin{array}{r}c_\rho s_\rho s_{\theta_m}(1 + c_{\Delta_m}) \\ + s_\rho c_{\theta_m}s_{\Delta_m}\end{array}
\end{vmatrix}
    \cr
& = s_\rho^2 s_{\Delta_m} = s_\rho^2\sin(\theta_{m + 1} - \theta_m)
    \cr
& = s_\rho^2\Im(e^{i\theta_{m + 1}}e^{-i\theta_{m}})
    \cr
& = s_\rho^2\Im\left[\frac{G(-p, m + 1, q)}{\sqrt q}\frac{\bar G(-p, m, q)}{\sqrt q}\right],
\end{align}

\noindent where we have used \eqref{e:G-pmq2} in the last line. On the other hand, the generalized quadratic Gau{\ss} sums can be explicitly calculated (see for instance the Appendix of \cite{HV2013}):
\begin{align}
G(-p, m, q) & = \sum_{l = 0}^{q - 1} e^{-2\pi i (p/q)l^2 + 2\pi i (m/q)l}
    \cr
& =
\begin{cases}
\sqrt q\jacobi{p}{q}e^{2\pi i \phi(p)m^2/q}, & \mbox{if } q \equiv 1 \bmod 4,
    \\
-i\sqrt q\jacobi{p}{q}e^{2\pi i \phi(p)m^2/q}, & \mbox{if } q \equiv 3 \bmod 4,
\end{cases}
\end{align}

\noindent where $\phi(p)$ denotes the inverse of $4p$ in the finite ring $\mathbb Z_q = \{0, 1, \ldots, q - 1\}$. Bearing in mind \eqref{e:G-pmq2}, and that the Jacobi symbol satisfies $\jacobi{p}{q}^2 = 1$, we get
\begin{align}
\label{e:triple1}
\begin{vmatrix}
\T(\tfrac{2\pi m}{Mq}^-)
    \cr
\T(\tfrac{2\pi m}{Mq}^+)
    \cr
\T(\tfrac{2\pi (m+1)}{Mq}^+)
\end{vmatrix}
& = s_\rho^2\Im\left[e^{2\pi i \phi(p)(m+1)^2/q}e^{-2\pi i \phi(p)m^2/q}\right]
    \cr
& = s_\rho^2\sin\left(\frac{2\pi \phi(p)(2m+1)}{q}\right),
\end{align}

\noindent where $s_\rho^2 = 1 - c_\rho^2$ is obtained from \eqref{e:cosrho}. The other quantity we are interested in is the scalar product of $\T(\tfrac{2\pi m}{Mq}^-) = (1, 0, 0)$ and $\T(\tfrac{2\pi (m + 1)}{Mq}^+)$:
\begin{align}
\label{e:scalar1}
\T(\tfrac{2\pi m}{Mq}^-)\cdot\T(\tfrac{2\pi (m + 1)}{Mq}^+) &  = c_\rho^2 -s_\rho^2c_{\Delta_m}
    \cr
& = c_\rho^2 -s_\rho^2\cos(\theta_{m + 1} - \theta_m)
    \cr
& = c_\rho^2 -s_\rho^2\Re(e^{i\theta_{m + 1}}e^{-i\theta_{m}})
    \cr
& = c_\rho^2 -s_\rho^2\Re\left[\frac{G(-p, m + 1, q)}{\sqrt q}\frac{\bar G(-p, m, q)}{\sqrt q}\right]
    \cr
& = c_\rho^2 -s_\rho^2\Re\left[e^{2\pi i \phi(p)(m+1)^2/q}e^{-2\pi i \phi(p)m^2/q}\right]
    \cr
& = c_\rho^2 -s_\rho^2\cos\left(\frac{2\pi \phi(p)(2m+1)}{q}\right).
\end{align}

\noindent Finally, taking \eqref{e:triple1} and \eqref{e:scalar1} respectively as the real and imaginary parts of a complex number, we define
\begin{align}
\label{e:zqmp1}
z_{q,m}(p) & \equiv
\begin{vmatrix}
\T(\tfrac{2\pi m}{Mq}^-)
    \cr
\T(\tfrac{2\pi m}{Mq}^+)
    \cr
\T(\tfrac{2\pi (m+1)}{Mq}^+)
\end{vmatrix}
+ i\T(\tfrac{2\pi m}{Mq}^-)\cdot\T(\tfrac{2\pi (m + 1)}{Mq}^+)
    \cr
& = i\,c_\rho^2 - i\,s_\rho^2\exp\left(\frac{2\pi i \phi(p)(2m+1)}{q}\right).
\end{align}

\noindent Summarizing, fixed $q$ and $m$, $[\T(\tfrac{2\pi m}{Mq}^-), \T(\tfrac{2\pi m}{Mq}^+), \T(\tfrac{2\pi (m+1)}{Mq}^+)]$, $\T(\tfrac{2\pi m}{Mq}^-)\cdot\T(\tfrac{2\pi (m + 1)}{Mq}^+)$, and, hence, $z_{q,m}(p)$,  depend exclusively on $\phi(p)$, i.e., on the inverse of $4p$ modulo $q$, which, as we will see in Section \ref{s:randomness}, is essentially random.

\subsection{Case with $q\equiv2\bmod4$}

\label{s:q24}

The cases with $q$ even are similar, so we will omit most details. When $q$ is even and $q/2$ is odd, $\psi(s, t_{pq})$ in \eqref{psistpq} adopts over the first period the form
\begin{equation}
\psi(s, t_{pq}) = \rho\sum_{m = 0}^{q/2 - 1}e^{i \theta_{2m + 1}} \delta(s - \tfrac{4\pi m + 2\pi}{Mq}), \quad s\in[0, \tfrac{2\pi}{M}),
\end{equation}

\noindent i.e., only the odd addends are to be considered. In this case, the vertices of $\X$, denoted by $\X_{2m+1}$, are located at $s = \tfrac{4\pi m + 2\pi}{Mq}$; so we have to calculate $[\T(\tfrac{2\pi (2m-1)}{Mq}^-), \T(\tfrac{2\pi (2m-1)}{Mq}^+), \T(\tfrac{2\pi (2m+1)}{Mq}^+)]$ and $\T(\tfrac{2\pi (2m-1)}{Mq}^-)\cdot\T(\tfrac{2\pi (2m+1)}{Mq}^+)$. The case $t = t_{pq} = t_{12}$ is trivial, with the first quantity being zero, and the second one being $\cos(\frac{4\pi}{M})$; hence, we consider $q > 2$.

The calculations for the triple product are exactly the same as in \eqref{e:triple1}, bearing in mind that we have to consider the right subscripts, i.e., substitute $c_m$ and $s_m$ by $c_{2m-1}$ and $s_{2m-1}$, respectively, and redefine $\Delta_m = \theta_{2m+1}-\theta_{2m-1}$. Therefore,
\begin{align}
\label{e:triple2a}
\begin{vmatrix}
\T(\tfrac{2\pi (2m-1)}{Mq}^-)
    \cr
\T(\tfrac{2\pi (2m-1)}{Mq}^+)
    \cr
\T(\tfrac{2\pi (2m+1)}{Mq}^+)
\end{vmatrix}
& = s_\rho^2 s_{\Delta_m}
    \cr
& = s_\rho^2\sin(\theta_{2m + 1} - \theta_{2m - 1})
    \cr
& = s_\rho^2\Im(e^{i\theta_{2m + 1}}e^{-i\theta_{2m - 1}})
    \cr
& = s_\rho^2\Im\left[\frac{G(-p, 2m + 1, q)}{\sqrt {2q}}\frac{\bar G(-p, 2m - 1, q)}{\sqrt {2q}}\right].
\end{align}

\noindent The generalized quadratic Gau{\ss} sums are now given by
\begin{align}
G(-p, 2m+1, q) & = 2G(-2p, 2m+1, q/2)
    \cr
& =
\begin{cases}
\sqrt{2q}\jacobi{2p}{q/2}e^{4\pi i \phi_1(p)(2m+1)^2/q}, & \mbox{if } q \equiv 2 \bmod 8,
    \\
-i\sqrt{2q}\jacobi{2p}{q/2}e^{4\pi i \phi_1(p)(2m+1)^2/q}, & \mbox{if } q \equiv 6 \bmod 8,
\end{cases}
\end{align}

\noindent where $\phi_1(p)$ is the inverse of $4(2p)=8p$ in $\mathbb Z_{q/2}$. Bearing in mind that $4(2m+1)^2-4(2m-1)^2 = 32m$, \eqref{e:triple2a} becomes
\begin{equation}
\begin{vmatrix}
\T(\tfrac{2\pi (2m-1)}{Mq}^-)
    \cr
\T(\tfrac{2\pi (2m-1)}{Mq}^+)
    \cr
\T(\tfrac{2\pi (2m+1)}{Mq}^+)
\end{vmatrix}
 = s_\rho^2\sin\left(\frac{32\pi \phi_1(p)m}{q}\right),
\end{equation}

\noindent where $s_\rho^2 = 1 - c_\rho^2$ is obtained from \eqref{e:cosrho}. On the other hand, $(8p)\phi_1(p) \equiv 1 \bmod (q/2)$ implies that $8\phi_1(p)$ is the inverse of $p$ in $\mathbb Z_{q/2}$, which we denote by $\phi(p)$. Therefore,
\begin{equation}
\label{e:triple2}
\begin{vmatrix}
\T(\tfrac{2\pi (2m-1)}{Mq}^-)
    \cr
\T(\tfrac{2\pi (2m-1)}{Mq}^+)
    \cr
\T(\tfrac{2\pi (2m+1)}{Mq}^+)
\end{vmatrix}
 = s_\rho^2\sin\left(\frac{2\pi \phi(p)m}{q/2}\right),
\end{equation}

\noindent where we prefer to write $q/2$ in the denominator, because we are working in $\mathbb Z_{q/2}$. Reasoning in the same way, the equivalent of \eqref{e:scalar1} is
\begin{equation}
\label{e:scalar2}
\T(\tfrac{2\pi (2m - 1)}{Mq}^-)\cdot\T(\tfrac{2\pi (2m + 1)}{Mq}^+) = c_\rho^2 + (c_\rho^2 - 1)\cos\left(\frac{2\pi \phi(p)m}{q/2}\right),
\end{equation}

\noindent and of \eqref{e:zqmp1} is
\begin{align}
\label{e:zqmp2}
z_{q,m}(p) & \equiv
\begin{vmatrix}
\T(\tfrac{2\pi (2m-1)}{Mq}^-)
    \cr
\T(\tfrac{2\pi (2m-1)}{Mq}^+)
    \cr
\T(\tfrac{2\pi (2m+1)}{Mq}^+)
\end{vmatrix}
+ i\T(\tfrac{2\pi (2m - 1)}{Mq}^-)\cdot\T(\tfrac{2\pi (2m + 1)}{Mq}^+)
    \cr
& = i\,c_\rho^2 - i\,s_\rho^2\exp\left(\frac{2\pi i \phi(p)m}{q/2}\right).
\end{align}

\noindent Summarizing, all the quantities depend exclusively on $\phi(p)$, i.e., on the inverse of $p$ modulo $q/2$.

\subsection{Case with $q\equiv0\bmod4$}

\label{s:q04}

When $q/2$ is even, $\psi(s, t_{pq})$ in \eqref{psistpq} adopts over the first period the form
\begin{equation}
\psi(s, t_{pq}) = \rho\sum_{m = 0}^{q/2 - 1}e^{i \theta_{2m}} \delta(s - \tfrac{4\pi m}{Mq}), \quad s\in[0, \tfrac{2\pi}{M}),
\end{equation}

\noindent i.e., only the even addends are to be considered. In this case, the vertices of $\X$, denoted by $\X_{2m}$, are located at $s = \tfrac{2\pi (2m)}{Mq}$; so we have to calculate $[\T(\tfrac{2\pi (2m)}{Mq}^-), \T(\tfrac{2\pi (2m)}{Mq}^+), \T(\tfrac{2\pi (2m+2)}{Mq}^+)]$ and $\T(\tfrac{2\pi (2m)}{Mq}^-)\cdot\T(\tfrac{2\pi (2m+2)}{Mq}^+)$. For that, we have to substitute in \eqref{e:triple1} $c_m$ and $s_m$ by $c_{2m}$ and $s_{2m}$, respectively, and redefine $\Delta_m = \theta_{2m+2}-\theta_{2m}$. Therefore,
\begin{align}
\label{e:triple3a}
\begin{vmatrix}
\T(\tfrac{2\pi (2m)}{Mq}^-)
    \cr
\T(\tfrac{2\pi (2m)}{Mq}^+)
    \cr
\T(\tfrac{2\pi (2m+2)}{Mq}^+)
\end{vmatrix}
& = s_\rho^2 s_{\Delta_m} = s_\rho^2\sin(\theta_{2m + 2} - \theta_{2m})
    \cr
& = s_\rho^2\Im(e^{i\theta_{2m + 2}}e^{-i\theta_{2m}})
    \cr
& = s_\rho^2\Im\left[\frac{G(-p, 2m + 2, q)}{\sqrt {2q}}\frac{\bar G(-p, 2m, q)}{\sqrt {2q}}\right].
\end{align}

\noindent In this occasion, the generalized quadratic Gau{\ss} sums are slightly more involved. Let us decompose $q = 2^rq'$, where $q'$ is odd; then,
\begin{equation}
G(-p, 2m, q) = G(-p, 2m, 2^rq') = G(-2^rp, 2m, q')G(-q'p, 2m, 2^r).
\end{equation}

\noindent On the one hand,
\begin{equation}
G(-2^rp, 2m, q') =
\begin{cases}
\sqrt {q'}\jacobi{2^rp}{q'}e^{2\pi i \phi_1(p)(2m)^2/q'}, & \mbox{if } q' \equiv 1 \bmod 4,
    \\
-i\sqrt {q'}\jacobi{2^rp}{q'}e^{2\pi i \phi_1(p)(2m)^2/q'}, & \mbox{if } q' \equiv 3 \bmod 4,
\end{cases}
\end{equation}

\noindent where $\phi_1(p)$ is the inverse of $4(2^rp)=2^{r+2}p$ in $\mathbb Z_{q'}$. On the other hand,
\begin{equation}
G(-q'p, 2m, 2^r) =
e^{\pi i \phi_2(p)(2m)^2/2^{r+1}}\jacobi{2^r}{q'p}(1 - i^{q'p})\sqrt {2^r},
\end{equation}

\noindent where $\phi_2(p)$ is the inverse of $q'p$ in $\mathbb Z_{2^r}$. Putting all together,
\begin{align}
\label{e:GGbar}
\frac{G(-p, 2m + 2, q)}{\sqrt {2q}}\frac{\bar G(-p, 2m, q)}{\sqrt {2q}} & = e^{2\pi i \phi_1(p)(2m+2)^2/q'}e^{\pi i \phi_2(p)(2m + 2)^2/2^{r+1}}
    \cr
& \hphantom{{}={}} \cdot e^{-2\pi i \phi_1(p)(2m)^2/q'}e^{-\pi i \phi_2(p)(2m)^2/2^{r+1}}
    \cr
& = e^{2\pi i [2^{r+2}\phi_1(p) + q'\phi_2(p)](2m+1) / q}.
\end{align}

\noindent This last expression can be further simplified. Indeed,
\begin{equation}
\begin{split}
2^{r+2}p\,\phi_1(p) \equiv 1 \bmod q' & \Rightarrow 2^{r+2}p\,\phi_1(p) + q'p\,\phi_2(p) \equiv 1 \bmod q'
    \\
q'p\,\phi_2(p) \equiv 1 \bmod 2^r & \Rightarrow 2^{r+2}p\,\phi_1(p) + q'p\,\phi_2(p) \equiv 1 \bmod 2^r;
\end{split}
\end{equation}

\noindent then, applying the well-known Chinese remainder theorem,
\begin{equation}
2^{r+2}p\,\phi_1(p) + q'p\,\phi_2(p) \equiv p(2^{r+2}\phi_1(p) + q'\phi_2(p)) \equiv 1 \bmod q,
\end{equation}

\noindent if and only if $2^{r+2}\phi_1(p) + q'\phi_2(p)$ is the inverse of $p$ modulo $q$, which we denote by $\phi(p)$:
\begin{equation}
2^{r+2}\phi_1(p) + q'\phi_2(p)\equiv\phi(p)\bmod q.
\end{equation}

\noindent Inserting this last expression into \eqref{e:GGbar}, it follows from \eqref{e:triple3a} that
\begin{equation}
\label{e:triple3}
\begin{vmatrix}
\T(\tfrac{2\pi (2m)}{Mq}^-)
    \cr
\T(\tfrac{2\pi (2m)}{Mq}^+)
    \cr
\T(\tfrac{2\pi (2m+2)}{Mq}^+)
\end{vmatrix}
= s_\rho^2\sin\left(\frac{2\pi\phi(p)(2m+1)}{q}\right),
\end{equation}

\noindent where $s_\rho^2 = 1 - c_\rho^2$ is obtained from \eqref{e:cosrho}. In the same way, the equivalent of \eqref{e:scalar1} is
\begin{equation}
\label{e:scalar3}
\T(\tfrac{2\pi (2m)}{Mq}^-)\cdot \T(\tfrac{2\pi (2m+2)}{Mq}^+) = c_\rho^2 + (c_\rho^2 - 1)\cos\left(\frac{2\pi\phi(p)(2m+1)}{q}\right),
\end{equation}

\noindent and of \eqref{e:zqmp1} is
\begin{align}
\label{e:zqmp3}
z_{q,m}(p) & \equiv
\begin{vmatrix}
\T(\tfrac{2\pi (2m)}{Mq}^-)
    \cr
\T(\tfrac{2\pi (2m)}{Mq}^+)
    \cr
\T(\tfrac{2\pi (2m+2)}{Mq}^+)
\end{vmatrix}
+ i\T(\tfrac{2\pi (2m)}{Mq}^-)\cdot \T(\tfrac{2\pi (2m+2)}{Mq}^+)
    \cr
& = i\,c_\rho^2 - i\,s_\rho^2\exp\left(\frac{2\pi i\phi(p)(2m+1)}{q}\right).
\end{align}

\noindent Summarizing, all the quantities depend exclusively on $\phi(p)$, i.e., on the inverse of $p$ modulo $q$.

We can combine the results of Sections \ref{s:q12}, \ref{s:q24} and \ref{s:q04} into the following theorem:

\begin{theorem}

\label{t:teo}

Let us consider the triple product of three consecutive tangent vectors (given by \eqref{e:triple1}, \eqref{e:triple2} and \eqref{e:triple3}), and the scalar product of a tangent vector and the second next one (given by \eqref{e:scalar1} and \eqref{e:scalar2} and \eqref{e:scalar3}). Then, those quantities depend exclusively on $\phi(p)$:
\begin{equation}
\label{e:phip}
\phi(p) \equiv
\begin{cases}
(4p)^{-1}\bmod q, & \mbox{if $q\equiv1\bmod 2$},
    \\
p^{-1}\bmod (q/2), & \mbox{if $q\equiv2\bmod 4$},
    \\
p^{-1}\bmod q, & \mbox{if $q\equiv0\bmod 4$}.
\end{cases}
\end{equation}

\noindent Furthermore, taking the first quantity as the real part and the second quantity as the imaginary part of a complex number $z_{q, m}(p)$ (defined respectively in \eqref{e:zqmp1}, \eqref{e:zqmp2} and \eqref{e:zqmp3}), $z_{q, m}(p)$ lies, for all $p$, on a circumference of center $ic_\rho^2$ and radius $s_\rho^2$, where $c_\rho$ is given by \eqref{e:cosrho}:
\begin{equation}
\label{e:zqmp}
z_{q, m}(p) =
\begin{cases}
i\,c_\rho^2 - i\,s_\rho^2\exp\left(\frac{2\pi i\phi(p)(2m+1)}{q}\right), & \mbox{if $q\not\equiv2\bmod4$, }
    \\
i\,c_\rho^2 - i\,s_\rho^2\exp\left(\frac{2\pi i\phi(p) m}{q/2}\right), & \mbox{if $q\equiv2\bmod4$. }
\end{cases}
\end{equation}

\end{theorem}

\section{Randomness}

\label{s:randomness}

In the previous section, we have shown in Theorem \ref{t:teo} how the triple product of three consecutive tangent vectors, and the scalar product of a tangent vector and the second next one, depend exclusively on $\phi(p)$ as defined in \eqref{e:phip}, i.e., we have  transformed the problem into one of finding inverses in finite rings. The existence and uniqueness of $\phi(p)$ is guaranteed since $\gcd(p, q) = 1$. Hence, for any given $q$, both $p$ and $\phi(p)$ can take $\varphi(q)$ different values in the corresponding finite ring, where $\varphi(q)$ is Euler's totient function, which gives the amount of positive integers less than or equal to $q$ that are coprime to $q$ (in fact, when $q\equiv2\bmod4$, we have $\varphi(q) = \varphi(2(q/2) = \varphi(2)\varphi(q/2) = \varphi(q/2)$). On the other hand, $z_{q, m}(p)$ can take $\varphi(q / \gcd(q, 2m + 1))$ different values, if $q\not\equiv2\bmod4$; and $\varphi((q/2) / \gcd(q/2, m))$ different values, if $q\equiv2\bmod4$. Therefore, we are interested in choosing $m$ such that $z_{q, m}(p)$ gives the largest possible amount of different numbers, i.e., such that $\gcd(q, 2m + 1) = 1$, if $q\not\equiv2\bmod4$; and $\gcd(q/2,m)=1$, if $q\equiv2\bmod4$. Without loss of generality, we can take $m = 0$, if $q\not\equiv2\bmod4$; and $m = 1$, if $q\equiv2\bmod4$. Then, \eqref{e:zqmp} becomes
\begin{equation}
\label{e:zqp}
z_q(p) =
\begin{cases}
i\,c_\rho^2 - i\,s_\rho^2\exp\left(\frac{2\pi i\phi(p)}{q}\right), & \mbox{if $q\not\equiv2\bmod4$, }
    \\
i\,c_\rho^2 - i\,s_\rho^2\exp\left(\frac{2\pi i\phi(p)}{q/2}\right), & \mbox{if $q\equiv2\bmod4$, }
\end{cases}
\end{equation}

\noindent which yields exactly $\varphi(q)$ different complex numbers lying on the same circumference. In other words, there is a one-to-one correspondence between $z_q(p)$ and $\phi(p)$. $\phi(p)$ can be efficiently computed, for instance, by the extended Euclidean algorithm; another more explicit (but less efficient) way is via Euler's theorem. For example, in $\mathbb Z_q$,
\begin{equation}
p^{\varphi(q)} \equiv 1\bmod q \Leftrightarrow p^{\varphi(q)-1} \equiv p^{-1}\bmod q, \quad\forall p\in\mathbb Z_q\slash\gcd(p,q) = 1.
\end{equation}

\noindent When $q$ prime, this last expression is known as Fermat's little theorem (of which Euler's theorem is in fact a generalization); in that case, $\varphi(q) = q - 1$, so
\begin{equation}
\label{e:fermat}
p^{q-1} \equiv 1\bmod q \Leftrightarrow p^{q-2} \equiv p^{-1}\bmod q, \quad\forall p\in\mathbb Z_q\backslash\{0\}.
\end{equation}

\noindent In this paper, however, we are not interested so much in finding $\phi(p)$, but rather in its randomness properties or, equivalently, in the randomness properties of \eqref{e:zqp}. There are diverse methods of generating pseudorandom numbers, the most popular ones being the linear congruential generators (LCGs) (see for instance \cite[Section 3.2.1]{knuth}). Given a large $q\in\mathbb N$ and $a, b, x_0\in\mathbb Z$, a linear congruential sequence $(x_n)_{n\ge0}$ of nonnegative integers smaller than $m$ is defined by
\begin{equation}
x_{n + 1} \equiv ax_n + b \bmod q, \quad n\ge 0.
\end{equation}

\noindent Then, after a careful choice of $q, a, b, x_0$, a sequence $(u_n)_{n\ge0}$ of linear congruential pseudorandom numbers uniformly distributed in the interval $[0, 1)$ is obtained by the normalization $u_n = x_n / q$, for $n\ge0$.

The quality of LCGs heavily depends on the coarseness of the lattice structure of $s$-dimensional vectors $\mathbf u_n^{(s)} = (u_n, \ldots, u_{n + s - 1})$ generated from the periodic sequence $(u_n)_{n\ge 0}$. It often happens \cite{entacher1998} that the lattice can be covered by a small amount of parallel hyperplanes: a sadly well-known example is the formerly popular RANDU generator
\begin{equation}
x_{n + 1} \equiv 65539 x_n \bmod 2^{31}.
\end{equation}

\noindent Since RANDU satisfies $x_{n + 2} \equiv 6 x_{n + 1} - 9 x_n \bmod 2^{31}$, it fails most three-dimensional criteria for randomness. Indeed, taking $(x_n, x_{n + 1}, x_{n + 2})$ as ``random'' points in the three-dimensional space, these points lie in exactly 15 planes! Therefore, the results obtained through RANDU are to be seen as suspicious.

In order to solve de deficiencies of LCGs, nonlinear random generators have been introduced \cite{eichenauer1988}. Their idea is that, given a large $q$ prime number, the elements are generated recursively by means of an integer-valued nonlinear function $f$:
\begin{equation}
x_{n + 1} \equiv f(x_n) \bmod q, \quad n\ge 0;
\end{equation}

\noindent then, we apply again the normalization $u_n = x_n/q$ as in the LCGs, to obtain pseudorandom numbers uniformly distributed over $[0, 1)$. An important particular case are the inversive congruential generators (ICGs), introduced by \cite{eichenauer1986}:
\begin{equation}
\label{e:ICG}
x_{n + 1} \equiv
\begin{cases}
a\, x_n^{-1} + b \bmod q, & x_n \ge 1,
    \\
b, & x_n = 0,
\end{cases}
\quad n\ge 0,
\end{equation}

\noindent with $q$ prime, $a\not\equiv0\bmod q$, which are characterized by the absence of any lattice structure, although their computational generation is not so efficient as with the LCGs. Remark that, in the literature, it is customary to write
\begin{equation}
x_{n + 1} \equiv a \overline x_n + b \bmod q, \quad n \ge 0,
\end{equation}

\noindent where $\overline z \equiv z^{p - 2} \bmod q$. From \eqref{e:fermat}, $\overline z$ is simply the multiplicative inverse of $z$, if $z\not\equiv0\bmod q$; while $\overline z$ is zero, if $z\equiv0\bmod q$.

Due to Eichenauer-Herrmann \cite{eichenauer1993} are as well the related explicit inversive congruential generators (EICGs), which are the relevant ones in this paper:
\begin{equation}
\label{e:EICG}
x_n \equiv \overline{an + b} \bmod q, \quad n \ge 0,
\end{equation}

\noindent with $q$ prime, $a\not\equiv0\bmod q$. It is immediate to see that $x_n$ has a period equal to $q$, i.e., $\{x_0, \ldots, x_{q-1}\} = \mathbb Z_q$; hence, any EICG with the normalization $u_n = x_n / q$ passes the uniformity test for equidistribution in $[0, 1)$. However, statistical independence properties of pseudorandom numbers are as important for stochastic simulations as uniformity properties. To study their statistical independence, Eichenauer-Herrmann used in \cite{eichenauer1993} the so-called serial test, which analyzes the discrepancy of tuples of pseudorandom numbers, and which we explain briefly here. The idea is, for a given dimension $k\ge 2$ and for $N$ arbitrary points $(\xi_0, \ldots, \xi_{N-1})\in[0, 1)^k$, to consider their discrepancy, which is defined as
\begin{equation}
D_N(\xi_0, \ldots, \xi_{N-1}) = \sup_J|F_N(J) - V(J)|,
\end{equation}

\noindent where the supremum is extended over all the subintervals $J$ of $[0, 1)^k$; $F_N(J)$ is $N^{-1}$ times the number of terms among $\xi_0, \ldots, \xi_{N-1}$ falling into $J$; and $V(J)$ denotes the volume of $J$.

In \cite{eichenauer1993}, given a sequence of numbers $(u_n)_{n\ge0}$ obtained with an EICG, the $k$-dimensional points
\begin{equation}
\mathbf u_n = (u_{n+n_1}, \ldots, u_{n+n_k})\in[0,1)^k, \quad 0\le n<p,
\end{equation}

\noindent were considered, with $n_1, \ldots, n_k$ arbitrary integers satisfying $0 = n_1 < \ldots < n_k < p$, and the abbreviation
\begin{equation}
D_p^{(k)} = D_p(\mathbf u_0, \ldots, \mathbf u_{p-1})
\end{equation}

\noindent being used for the discrepancy of the points. Then, an EICG passes the $k$-dimensional serial test if $D_p^{(k)}$ is reasonably small. In this regard, the following two theorems were formulated in \cite{eichenauer1993}:
\begin{theorem}
\label{t:t1}
Let $2\le k< p$. Then, the discrepancy $D_p^{(k)}$ for any EICG satisfies
$$
D_p^{(k)} < 2p^{-1/2}\left((k-1)\left(\frac{2}{\pi}\log p + \frac{7}{5}\right)^k+1\right)+kp^{-1}.
$$
\end{theorem}

\begin{theorem}
\label{t:t2}
Let $0 < t \le 1$. Then there exist more than $A_p(t)(p-1)$ values of $a\in\mathbb Z_p^*$ such that the discrepancy $D_p^{(k)}$ for any corresponding EICG satisfies
$$
D_p^{(k)}\ge\frac{t}{2(\pi + 2)}p^{-1/2}
$$

\noindent for all dimensions $k \ge 2$, where
$$
A_p(t) = \frac{(1-t^2)p}{(4-t^2)p + 12p^{1/2} + 9}.
$$

\end{theorem}

Theorems \ref{t:t1} and \ref{t:t2} show that, in the EICG method, the discrepancy $D_p^{(k)}$ has on the average an order of magnitude between $p^{-1/2}$ and $p^{-1/2}(\log p)^k$. However, it is precisely in this range of magnitudes where the discrepancy of $p$ independent and uniformly distributed points taken from $[0, 1)^k$ is found, which is roughly $p^{-1/2}(\log\log p)^{1/2}$. In this sense, we can say that EICGs model true random numbers very closely, or, in Eichenauer-Herrmann's words, EICGs have \emph{even better structural and statistical independence properties than the standard type}, i.e., than \eqref{e:ICG}. Furthermore, they also behave very well in parallel and vector computations, as shown by Niederreiter in \cite{niederreiter}. Indeed, if we define a family of $N$ EICGs:
\begin{equation}
x_n^i \equiv \overline{a^in + b^i} \bmod q, \quad n \ge 0, \quad  i = 1, \ldots, N,
\end{equation}

\noindent then, the $N$-tuples of the form $(x_n^1, \ldots x_n^N)$ have good statistical properties if all the $N$ numbers $b^i\overline{a^i}$ are distinct. Summarizing, this approach is, in Niederreiter's words, \emph{eminently suitable for the generation of parallel streams of pseudorandom numbers with desirable properties}.

In principle, it could be possible to work in $\mathbb Z_q$, with $q$ an arbitrary natural number, although the most common choices are $q$ prime, as in the definitions \eqref{e:ICG} and \eqref{e:EICG}, or $q$ a power of two. For instance, in \cite{eichenauer1994}, Eichenauer and Ickstadt study the equally good pseudorandom properties of EICGs defined by the inverses of the odd integers in $\mathbb Z_q$, with $q = 2^\omega$, $\omega\ge 5$:
\begin{equation}
\label{e:EICGeven}
x_n \equiv (an + b)^{-1}\bmod 2^\omega, \quad n \ge 0,
\end{equation}

\noindent with $a\equiv2\bmod4$, $b\equiv1\bmod2$; it is immediate to see that $x_n$ has a period equal to $q/2=2^{\omega-1}$, i.e., $x_n$ takes all the possible odd values in $\mathbb Z_q$, or, in other words, $\{x_0, \ldots, x_{q/2-1}\} = \mathbb Z_q^*$.

Coming back to \eqref{e:zqp}, all the previous arguments should be more than enough to justify the extremely good pseudorandom character of \eqref{e:zqp} and, hence, of $\X$ and $\T$. In particular, when $q$ is an odd prime,
\begin{equation}
\label{e:zqpprime}
z_q(p) = i\,c_\rho^2 - i\,s_\rho^2\exp(2\pi i u_p),
\end{equation}

\noindent where $u_p = x_p / q$, and $x_p$ is given by \eqref{e:EICG}, with $a \equiv 4\bmod q$, $b \equiv 0\bmod q$; i.e., $x_p\equiv \overline{4p}\bmod q$. Therefore, by direct application of Theorems \ref{t:t1} and \ref{t:t2}, $u_p$ is a sequence of pseudorandom numbers uniformly distributed in the interval $(0, 1)$, which, from our Theorem \ref{t:teo}, implies that $z_q(p)$ is a sequence of pseudorandom numbers uniformly distributed in the circumference of center $i\,c_\rho^2$ and radius $s_\rho^2$. Observe that we have to exclude the case $p = 0$ and, hence, the numbers lie on $(0, 1)$ instead of $[0, 1)$. Nevertheless, since $u_0 = 0$, we are just omitting the first term of the sequence (and the point $z_0 = i\cos(2\rho)$ in the circumference), which makes this minor issue irrelevant for all purposes.

The same reasoning is valid when $q$ is twice a prime number. Then, \eqref{e:zqpprime} also holds, with $u_p = x_p / (q / 2)$ and $x_p\equiv \overline{p}\bmod (q/2)$, i.e., we are taking $a\equiv 1\bmod (q/2)$ and $b\equiv 0\bmod (q/2)$ in \eqref{e:EICG}. Again, we exclude $u_0 = 0$, and the whole previous paragraph is valid in its integrity.

Let us mention also the case with $q = 2^\omega$. Then, \eqref{e:zqpprime} also holds, with $u_p = x_p / q$, and $x_p \equiv \overline{2p - 1} \bmod q$, i.e, we are taking $a\equiv 1\bmod q$ and $b\equiv-1\bmod q$ in \eqref{e:EICGeven}. Again, we have obtained a sequence of pseudorandom numbers uniformly distributed in the circumference of center $i\,c_\rho^2$ and radius $s_\rho^2$, but, unlike the two previous cases, it is not necessary to exclude any number.

Analyzing all the possible values of $q$ lies certainly beyond the scope of this paper. Moreover, \eqref{e:zqpprime} is not the only possible probability generator which can be obtained from the evolution of $\X$ in $\T$. For instance, given $N$ different primes $q_1, \ldots, q_N \ge 5$, it is possible to combine \eqref{e:zqpprime} in a way that closely resembles the so-called \emph{compound approach} explained in \cite{eichenauer1993b}. Let us particularize \eqref{e:zqp} as
\begin{equation}
\frac{c_{\rho_j}^2 + iz_{q_j}(p)}{s_{\rho_j}^2} = \exp\left(\frac{2\pi i\phi_j(p)}{q_j}\right),
\end{equation}

\noindent where $\rho_j$ is the angle corresponding to $q_j$; $\phi_j(p) \equiv (4p)^{-1}\bmod q_j$; and $p\not\equiv0\bmod q_j$. Then,
\begin{equation}
\label{e:prodzqp}
\prod_{j = 1}^N\frac{c_{\rho_j}^2 + iz_{q_j}(p_j)}{s_{\rho_j}^2} = \exp\left(2\pi i\sum_{j = 1}^N\frac{\phi_j(p)}{q_j}\right).
\end{equation}

\noindent Denoting now
\begin{equation}
u_p \equiv \sum_{j = 1}^N\frac{\phi_j(p)}{q_j}\bmod 1,
\end{equation}

\noindent \eqref{e:prodzqp} becomes
\begin{equation}
\prod_{j = 1}^N\frac{c_{\rho_j}^2 + iz_{q_j}(p_j)}{s_{\rho_j}^2} = \exp(2\pi iu_p), \qquad p\not\equiv0\bmod q_j, \forall j.
\end{equation}

\noindent The left-hand side is directly obtained from $\T$, and its good random properties follows directly from \cite{eichenauer1993b}. The only minor difference is that, in our case, $p$ can take $(q_1-1)\cdot\ldots\cdot(q_N-1)$ different values modulo $q_1\cdot\ldots\cdot q_N$, while its equivalent in \cite{eichenauer1993b} can take all the $q_1\cdot\ldots\cdot q_N$ values. However, in practice, taking $q_1,\ldots, q_N$ large enough, the amount of values that we are excluding is, for all purposes, irrelevant.

\section{Conclusions}

\label{s:conclusions}

In this paper, we have considered the evolution of \eqref{e:xt}-\eqref{e:schmap}, taking a regular planar polygon of $M$ sides as the initial datum. Bearing in mind the recent results in \cite{HV2013}, where we gave very strong evidence that $\X(s, t)$ is a skew polygonal at times which are rational multiples of $2\pi/M^2$; we have studied \eqref{e:xt}-\eqref{e:schmap} from a completely novel point of view: that of an evolution equation which yields a very good pseudorandom generator in a completely natural way.

Due to the algebraic complexity of the calculations involved, we have limited ourselves mainly to the study at rational times of two quantities, which are illustrative enough of the essential random character of \eqref{e:xt}-\eqref{e:schmap}: the triple product of three consecutive tangent vectors; and the scalar product of a tangent vector and the second next one. These quantities, when taken respectively as the real and imaginary parts of a complex number, yield an excellent generator of pseudorandom numbers uniformly located on a circumference. Furthermore, it is straightforward to combine different rational times to develop additional pseudorandom generators.

Although the main aim of this paper is to show the randomness in the evolution of \eqref{e:xt}-\eqref{e:schmap}, for which it is largely enough to work with $\T$, it is not irrelevant to mention here that, as observed in \cite{HV2013}, $\X(0, t)$ is very intimately related to Riemann's nondifferentiable function,
\begin{equation}
f(t) = \sum_{k = 1}^{\infty}\frac{\sin(\pi k^2 t)}{\pi k^2},
\end{equation}

\noindent which, as proved by Jaffard \cite{jaffard1996}, is a multifractal and, in fact, fits under the so called Frisch-Parisi conjecture proposed in \cite{frischparisi1985} (see also \cite{frisch1995} for more details). Therefore, although giving a complete algebraic characterization of $\X(0, t)$ reveals as a complex task which clearly lies beyond the scope of this paper and which we postpone for the future, we can nonetheless expect an even richer randomness structure in $\X$.

The ideas presented here can be most probably extended to other types of evolution equations, in order to obtain new probability generators. Obviously, this approach is not intended for competing with commercially developed algorithms; even though, during the simulation of \eqref{e:xt}-\eqref{e:schmap}, large sequences of pseudorandom numbers with good statistical properties can be generated with virtually no additional cost, i.e., for free.

As we have said in the introduction, a recurring question is up to what extent VFE is valid as a simplified model. In this line, the random character of \eqref{e:xt}-\eqref{e:schmap} proved in this paper is at the very least not in contradiction with the physical motion of a real vortex filament. Furthermore, we venture to suggest that finding the existence of well-behaved pseudorandom sequences of numbers \emph{inside} the evolution of a proposed physical model might be a first test in validating that model with respect to the phenomenon that it is trying to describe. Indeed, real natural phenomena are in general characterized by their chaotic, truly random behaviour. Therefore, a model with an easily predictable structure might be suspected not to match reality accurately.

\end{document}